\documentclass[letterpaper, 10 pt, conference]{ieeeconf}  
\IEEEoverridecommandlockouts                              
\usepackage{amssymb}  
\usepackage{amsmath}
\usepackage{array}
\usepackage[OT1]{fontenc} 
\usepackage[table,dvipsnames]{xcolor}
\usepackage{booktabs}
\usepackage{svg}
\usepackage{tikz}
\usepackage{pgfplots}
\pgfplotsset{compat=1.17}
\usepgfplotslibrary{groupplots}
\usetikzlibrary{spy,chains,shapes.multipart,shapes,calc,fit,automata,positioning}
\usepackage{acronym}
\usepackage[font=footnotesize]{caption}
\usepackage[font=footnotesize,labelformat=simple]{subcaption} 

\usepackage{mathtools}
\usepackage[intoc]{nomencl}
\usepackage{bbm,dsfont}
\DeclareMathOperator{\sgn}{sgn}
\DeclareMathOperator{\diag}{diag}
\DeclareMathOperator*{\argmin}{arg\,min}

\newcounter{mcntr}[table]
\newcounter{rowcntr}[table]
\renewcommand{\therowcntr}{\arabic{mcntr}\alph{rowcntr}}

\newcolumntype{N}{>{\refstepcounter{rowcntr}\therowcntr}c}

\newcommand{\calexp}[1]{\mathbb{E}\left[#1\right]}

\newtheorem{prob}{\textbf{Problem}}
\newtheorem{prop}{\textbf{Proposition}}
\newtheorem{rmk}{\textbf{Remark}}

\acrodef{ODE}       [ODE]           {Ordinary Differential Equation}
\acrodef{CTMC}	    [CTMC]	    {Continuous Time Markov Chain}
\acrodef{MDP}	    [MDP]	    {Markov Decision Process}
\acrodefplural{MDP}  [MDPs]	    {Markov Decision Processes}
\acrodef{PDE}	    [PDE]	    {Partial Differential Equation}
\acrodef{HJBE}	    [HJBE]          {Hamilton-Jacobi-Bellman Equation}
\acrodef{MJP}	    [MJP]          {Markov Jump Process}

\makenomenclature
\renewcommand{\nomname}{}

\title{\Large \bf Optimal control of stochastic networks of $M/M/\infty$ queues with linear costs}
\author{Giovanni Pugliese Carratelli and Ioannis Lestas
\thanks{\noindent G. Pugliese Carratelli and I. Lestas are with the Department of Engineering, University of Cambridge, United Kingdom.
 {\tt\small \{gp459, icl20\}@cam.ac.uk }}}%

\begin{document}
\maketitle
\thispagestyle{empty}
\pagestyle{empty}
\begin{abstract}\label{sec.Abstract}
We consider an arbitrary network of $M/M/\infty$ queues with controlled transitions between queues. We consider optimal control problems where the costs are linear functions of the state and inputs over a finite or infinite horizon. We provide in both cases an explicit characterization of the optimal control policies. We also show that these do not involve state feedback, but they depend on the network topology and system parameters. The results are also illustrated with various examples.
\end{abstract}
\section{Introduction}\label{sec.introduction}

Studies on the control of queues have appeared in the literature form an early stage \cite{Erlang1909}. Fundamental work on the control of arrival and service rates of single queues was established in studies such as \cite{Crabill1972,Stidham1985,Johansen1980} where the authors quantify the conditions for which the optimal control polices have a special monotone structure \cite{Topkis1978}. Networks of queues have also received significant attention due to the numerous applications in operations research \cite{Li2019}, communication networks as well as biological systems \cite{arazi2004}.
Network models \cite{Jackson1957} with multiple interacting queues have been studied for $M/M/1$ \cite{Ghoneim1985,Rosberg1982,Hajek1984} and $M/G/1$ queues \cite{BadianPessot2019,Gallisch1979, HernandezLerma1983} in specific configurations \cite{Weber1987}. More recently \cite{Su2024} examines the problem of double sided queues while \cite{Dimitrakopoulos2017}, \cite{Ata2006} and \cite{Adusumilli2010} considered service and arrival control problems in $M/M/1$ queues with various assumptions on costs in a \ac{MDP} setting.

The study of $M/M/\infty$ queues has received less attention with such examples including the work in \cite{Feinberg2015} where the authors focus upon cost related aspects for a parallel configuration.
In this study we characterise the optimal policies for stochastic networks of $M/M/\infty$ queues with linear costs, prescribed service rates, and routing events between queues with controlled rates. In particular, we derive explicit expressions for the optimal policies over both a finite and an infinite horizon. We also show that the optimal policies do not involve state feedback, but they depend on the topology of the network and the system parameters.

One of the critical aspects of queueing networks is the positivity of the systems states which we model via an appropriate \ac{MJP}.  In our analysis we derive  appropriate expressions for the \ac{HJBE} for the network and quantify explicitly the value function. Our work also has links to other types of problems \cite{Li2024,Blanchini2023} where results have been derived for deterministic positive linear systems with input signals that are constrained by a linear function of the state.

We would also like to note that the methodology used to derive the optimal policies is part of ongoing work with potential extensions to more broad classes of stochastic networks and corresponding optimal control problems.

The manuscript is organised as follows. In Section \ref{sec.Prelim} we introduce the notation and the models that will be considered. In Section \ref{sec.Probform} we define the optimal control problem we seek to address. Our main results are stated in Section \ref{sec.mainresults}. In Section \ref{sec.Methods} we provide the proofs of our results. In Sections \ref{sec.examples} and \ref{sec.Appendix} we provide examples validating our results.  Finally, the paper is concluded in Section \ref{sec.Conclusions}.

\section{System model and notation}\label{sec.Prelim}
In Section \ref{sec.Model} we present the mathematical model and quantities employed throughout the study.

\subsection{Notation}\label{ssec.Notation}
\printnomenclature
\markboth{\nomname}{\nomname}
\nomenclature{$\mathbb{R}_>$}{Set of positive real numbers $\{x\in\mathbb{R}:x >0\}$}
\nomenclature{$\mathbb{R}_\geq$}{Set of non-negative real numbers $\{x\in\mathbb{R}:x \geq0\}$}
\nomenclature{$\mathbb{Z}_>$}{Set of positive integers $\{1,2,3,\ldots\}$}
\nomenclature{$\mathbb{Z}_\geq$}{Set of non-negative integers $\{0,1,2,3,\ldots\}$}

\subsection{System model and mathematical preliminaries}\label{sec.Model}
We consider the \ac{MJP} representing a network of $n$ interacting queues.  The number of elements in each queue\footnote{For convenience in the presentation, we will often slightly abuse notation and refer to a queue $i$ as queue $X_i$.} $i$ at time $t$ is a random variable  $X_i(t)$ and we also denote $X(t)= [X_1(t),X_2(t),\ldots ,X_n(t)]$. Each element in a queue can undergo transitions which are one of the two types described in Table \ref{table:def}: exit/servicing, or transition to another queue. Each event occurs after an exponential time, with the rate for transition to other queues being controlled. The overall system is represented in \eqref{eq.events} below and is formally defined via its Kolmogorov equation, which is provided later in this Section in \eqref{eq.FKE}.
\begin{equation}\label{eq.events}
        X(t) \xrightarrow[]{W_i(X(t),t)} X(t)+r_i, \forall i \in \mathcal{I}=\{1,\ldots,m\}, m \in \mathbb{Z}_>
\end{equation}
In particular, we have $m$ possible discrete events. Each function $W_i: \mathbb{Z}_\geq^{n}\times \mathbb{R} \to \mathbb{R}_\geq$ denotes the rate of event $i$. The vector $r_i\in\mathbb{Z}^{n}$ denotes the change of the state of the system $X$ due to event $i$. This is also denoted as the $i-th$ column of a state change matrix $R\in \mathbb{Z}_\geq^{n\times m}$, \emph{i.e.} matrix $R=[r_i]_{i\in\mathcal{I}}$ is constructed by stacking side by side the column vectors $r_i$. It should be noted that in our system the events and their rates are such that the state $X(t)$ of the system remains non-negative.

The events we consider are summarised in Table \ref{table:def}. We consider $m_e\leq m, m_e\in \mathbb{Z}_\geq$ exit events from the network and $m_u\leq m, m_u\in \mathbb{Z}_\geq$ routing events that transition a unit from queue $X_i$ to queue $X_j$. The events correspond to the rows \ref{ev.Ex}, \ref{ev.Dist} in Table \ref{table:def} and completely describe the evolution of the system, \emph{i.e.} $m =m_u+m_e$. We partition the set of indices $\mathcal{I}$ associated with the considered events defined in \eqref{eq.events} in two mutually exclusive sets. The index set $\mathcal{E}, |\mathcal{E}|= m_e$ is associated with the $m_e$ exit events and the set $\mathcal{D}, |\mathcal{D}|=m_u$ is associated with the $m_u$ controlled routing events.
\begin{table}[!h]
       \rowcolors{2}                                 
           {black!10!white}    
           {}                                  
        \centering
        \begin{tabular}{Ncccc}
            \toprule
	    \multicolumn{1}{c}{\textbf{ID}}  &\textbf{Event} & \textbf{Transition} & \textbf{Rate} & \textbf{Index set}\\
                    \midrule
	    \label{ev.Dist}	&Routing& $(x_i,x_j) \xrightarrow[]{} (x_i-1,x_j+1)$ & $u_{i_j}x_i$ & $\mathcal{D}$\\
	    \label{ev.Ex}  	&Exit & \emph{$x_i \xrightarrow[]{} x_i-1$}  &  $\gamma_i x_i$& $\mathcal{E}$ \stepcounter{mcntr}\setcounter{rowcntr}{0}\\
                \bottomrule
        \end{tabular}
	\caption{Events associated with the considered system}\label{table:def}
\end{table}
The rate $W_i$ associated with the exit of units from a queue $i$ when $X_i(t)=x_i$  is $W_i(x,t)=\gamma_i x_i \ \forall i \in \mathcal{E}$ where $\gamma_i \in \mathbb{R}_\geq$, \emph{i.e.} each individual element in the queue has exit rate $\gamma_i$. For each routing event a single unit is transitioned from queue $i$ to queue $j$ $\forall i,j \in \mathcal{D}$ with $i\neq j$. The rate of this event at time $t$ when $X_i(t)=x_i$ is equal to the product of the queue size $x_i$ and the control variable $u_{i_j}(t)$; \emph{i.e} $W_{i_j}(x,t)= x_i u_{i_j}(t) \ \forall {i_j} \in \mathcal{D}$, corresponding to the fact that each individual element in the queue has transport rate $u_{i_j}$. Note that in the following Section we address optimal control problems where we optimise the control inputs over set $\mathcal{U}$ that is a constraint for the routing rates.

The system we consider is defined via the Kolmogorov Equation (also known as the Master equation). This is a partial difference equation for the probability at time $t$ the number of elements in each queue takes specific values. For any $x\in \mathbb{Z}_\geq^{n}$ we denote by $\mathbb{P}(x,t)$ the probability that $X(t) = x$. The master equation for the system is
\begin{align}\label{eq.FKE}
        \dfrac{\partial{\mathbb{P}(x,t)}}{\partial t}&=
 \sum_{j\in \mathcal{E}} W_j(x-r_j,t)\mathbb{P}(x-r_j,t)-W_j(x,t)\mathbb{P}(x,t)]\notag\\
& +\sum_{k\in \mathcal{D}} W_k(x-r_k,t)\mathbb{P}(x-r_k,t)-W_k(x,t)\mathbb{P}(x,t)]
 \end{align}

For convenience in the notation throughout the paper we also define the following quantities for each of the considered event types. We define matrices $R_\mathcal{D}\in\mathbb{Z}_\geq^{n\times m_u}$, $R_\mathcal{E}\in\mathbb{Z}_\geq^{n\times m_e}$, which are state vector change matrices for routing events and for the exit events respectively, as follows
\begin{equation}
	R_\mathcal{D} = [r_i]_{i \in\mathcal{D}},	R_\mathcal{E} = [r_i]_{i \in\mathcal{E}}
\end{equation}
where each column $r_i$ corresponds to the change of the state $X$ of the system for the corresponding event (as illustrated in \eqref{eq.events} and the text below it).
\nomenclature{$L = [l_i]_{i \in S}$}{For a set $S\subset \mathbb{Z}_>$ the matrix $L$ is defined by stacking side by side the column vectors $l_i \forall i \in S$}

The rate parameter matrices for the considered events are $\Gamma \in \mathbb{R}_>^{m_e\times m_e}$
and $U(t) \in \mathbb{R}_\geq^{m_u\times m_u}$ and are defined as
\begin{align}
&\Gamma = \diag[\gamma_1, \ldots, \gamma_{m_e}]\notag\\
& U(t) = \diag[u_1(t),\ldots, u_{m_u}(t)] \label{def:U}
\end{align}
\nomenclature{$A = \diag(v)$}{For vector $v$ of size $n$ then $A_{ij} = v_i \forall i,j \in 1,\ldots, n$ with $i=j$ and $A_{ij} =0 \forall i \neq j$}
These are diagonal matrices that include the routing and serving rate parameters, respectively.

We also define two one/zero matrices mapping the rate and control matrices to the specific queues in the network. In particular, we define the matrices $H\in \mathbb{Z}_\geq^{m_u\times n}$, $E\in \mathbb{Z}_\geq^{m_e\times n}$ where $H_{ki} =1$ if routing event $k$  has as source queue $i$ and $E_{ki}=1$ if exit event $k$ occurs at queue $i$.
\nomenclature{$\mathbb{E}[\cdot]$}{Expectation operator}

\section{Optimal control problem formulation and considered costs}\label{sec.Probform}
We consider the problem of finding an optimal feedback policy to minimise a trade-off between the cost for maintaining units in the network and control costs for a network of $M/M/\infty$ queues, as described in Section \ref{sec.Model}.
We seek to minimise the costs for the presence of $x$ units in the network and the costs for routing units in the network. We first consider the problem of minimising the total expected costs in continuous time over a finite time horizon $T\in \mathbb{R}_>$ and we then consider the problem of finding an optimal policy over an infinite horizon.

We take into consideration the stage cost $g_c(x(t),U(t))$ per unit of time
\begin{equation}\label{eq.stagecost}
	g_c(x(t),U(t)) = q^Tx(t) + v^TU(t)Hx(t)
\end{equation}
and we consider the following terminal cost 
\begin{equation}\label{eq.finalcost}
\tilde g_{c}(x(T)) = c^Tx(T)
\end{equation}
\nomenclature{$\cdot^T$}{Transpose of matrix or vector}
The constant $q \in \mathbb{R}_\geq^{n}$ denotes non-negative cost coefficients for the vector $x(t)$ of the number of units in the $n$ queues at time $t$ and $c \in \mathbb{R}_\geq^{n}$ are cost coefficients at the final time $T\in \mathbb{R}_>$ in the considered horizon. The vector $v \in \mathbb{R}_\geq^{m_u}$ denotes non-negative cost coefficients for the vector of control signals. The proportionality between $U$ and $x$ suggests a cost is incurred for each unit to maintain a prescribed routing rate. It should be noted that the costs of the system remain positive because the costs are the product of the non-negative state $x(t)$ of the system with non-negative costs and non-negative control signal.

We search over deterministic feedback policies that are a function of the current state of the system\footnote{We focus on searching over deterministic policies in this paper. The inclusion of randomized policies will be considered in future work.}. This is without loss of generality due to the Markov nature of the system, \emph{i.e.} state feedback policies would be optimal even if the policy was allowed to depend on the history of the process \cite{Blackwell1964,Puterman1994,Bertsekas2005}.  It should be noted that the set of policies we consider are constrained to take arbitrary non-negative values in $\mathcal{U}$ resulting in a constraint for the optimal control problems defined below.

We denote a particular state feedback policy for the system as in \eqref{def:U}, \emph{i.e.} $U(x,t)=\diag(u_k(x,t))$,  $u_k(x,t):\mathbb{Z}^{n}_\geq\times \mathbb{R}_\geq \to \mathcal{U}$  and we denote by $\mathcal{C}$ the set of all time varying state feedback control policies $U$ where $u_k$ take values in $\mathcal{U}$. We also denote by $\mathcal{\tilde C}$ the set of all time invariant state feedback control policies $U$.
We now consider the following finite horizon optimal control problem for the system described in Section \ref{sec.Model}.
\begin{prob}\label{prob.FH}
Consider the system in \eqref{eq.FKE} as described in Section \ref{sec.Model}, the stage cost \eqref{eq.stagecost} and the terminal cost \eqref{eq.finalcost}. We consider the following optimal control problem
\begin{equation}\label{eq.MinProbEqFH}
        \mathcal{V}(x_0,0) = \min_{U\in\mathcal{C}}\mathcal{V}_U(x_0,0)
\end{equation}
where $\mathcal{V}_U(x_0,0):\mathbb{Z}^n_\geq\times \mathbb{R}_\geq\to \mathbb{R}_\geq$ defined below is the cost of the evolution of the \ac{MJP} in \eqref{eq.FKE} from the initial condition~$x_0$.
\begin{equation}\label{eq.totcostFH}
        \mathcal{V}_U(x_0,0) = \calexp{c^Tx(T) + \int_0^Tg_c(x(t),U)dt\middle\vert X(0)=x_0}
\end{equation}
\end{prob}
The solution of Problem \ref{prob.FH} corresponds to finding the optimal policy $U^*(x,t)$ that minimises the total expected costs in \eqref{eq.totcostFH}. It should be noted that Problem \ref{prob.FH} is a finite horizon stochastic optimal control problem with constrained control variables and linear costs.

We also consider the infinite horizon optimal control problem stated below. Here we assume that  $X=0$ is an absorbing state of the system. That is for all queues $i$ there is a "directed path" to a queue $j$ with a non zero serving rate $\gamma_j \neq 0$, \emph{i.e.} there is a choice of control inputs such that there exists a non zero probability that a unit in queue $i$ can reach queue $j$ in finite time. It should be noted that this is a mild assumption associated with the well-posedness of the problem.
\begin{prob}\label{prob.IH}
Consider the system in \eqref{eq.FKE} as described in Section \ref{sec.Model}
and the stage cost \eqref{eq.stagecost}. We consider the following optimal control problem
\begin{equation}\label{eq.MinProbEqIH}
        \mathcal{V}(x_0) = \min_{U\in\mathcal{\tilde C}}\mathcal{V}_U(x_0)
\end{equation}
where $\mathcal{V}_U(x_0):\mathbb{Z}^n_\geq\to \mathbb{R}_\geq$ defined below is the cost of the evolution of the \ac{MJP} in \eqref{eq.FKE} from the initial condition $x_0$.
\begin{equation}\label{eq.totcostIH}
        \mathcal{V}_U(x_0) =\lim_{T\to +\infty}\calexp{  \int_0^Tg_c(x(t),U)dt\middle\vert X(0)=x_0}
\end{equation}
\end{prob}

\section{Results}\label{sec.mainresults}
We give our first result which provides an explicit characterization of the optimal policy for Problem \ref{prob.FH}. We also show that the optimal policy does not involve state feedback.
\begin{prop}\label{prop.FH}
Consider Problem \ref{prob.FH} for the system in \eqref{eq.FKE}  described in Section \ref{sec.Model} with stage costs \eqref{eq.stagecost}. Then the policy $U^*$ in \eqref{eq.optsol} is an optimal policy.
        \begin{equation}\label{eq.optsol}
                U^*(x,t) = \diag\left(\frac{u_{max}}{2}(\mathds{1}-\sgn(y^T(t){R_\mathcal{D}}+v^T))\right)
        \end{equation}
where $y^T(t)\in \mathbb{R}^{n}$ is a solution of
        \nomenclature{$\sgn(x)$}{$\sgn(x) = 1$ if $x>0$ and $\sgn(x) = -1$ if $x<0$. If $x$ is a vector then $\sgn(x)=[\sgn(x_1),\ldots, \sgn(x_n)]^T$}
\nomenclature{$\mathds{1}$}{Vector of ones}
\begin{align}\label{eq.paramState}
&\left[\dot y^T(t)+q^T + y^T(t)R_\mathcal{E}\Gamma E  \right.\\
&       \left.+ \frac{u_{max}}{2}\left((y^T(t)R_\mathcal{D}+v^T)-\lvert y^T(t)R_\mathcal{D}+v^T\rvert_{elem}\right)H\right]=0\notag
\end{align}
with the terminal condition
\begin{equation}
        y(T) = c
\end{equation}
\end{prop}
\begin{proof}
	See Section \ref{sec.Methods}.
\end{proof}

\begin{rmk}\label{rmk.indipend}
Proposition \ref{prop.FH} shows that the optimal policy $U^*$ is a time varying function that is independent of state $X(t)$ of the network. This implies there is no incentive to implement a state feedback scheme that dynamically adjusts the optimally chosen routing policy based on the current value of the state. We would also like to note that the optimal policy depends on the topology of the network and the system parameters through $y(t)$.
\end{rmk}
\begin{rmk}
It should be noted that Proposition \ref{prop.FH} holds for arbitrary system and cost parameters. In particular the result holds for any $\Gamma$ and $\mathcal{U}$ and for any non-negative cost coefficient vectors $q$ and $v$. Also our result holds for arbitrary network interconnection matrices $H$ and $E$ as long as \eqref{eq.paramState} has a solution.
\end{rmk}
\begin{rmk}
A direct corollary of Proposition \ref{prop.FH} is that the total cost for Problem \ref{prob.FH} is $\mathcal{V}(x_0,0) = y^T(0)x_0$.
\end{rmk}
\begin{rmk}
It should be noted that one can considerably reduce the computation time to obtain the optimal policy established in Proposition \ref{prop.FH} by direct integration of \eqref{eq.paramState}
\end{rmk}
\begin{rmk}
	Proposition \ref{prop.FH} can be extended to the case where the maximum control rate is linearly dependent with the state or varies with time and also when a different bound is used for each control input $u_k$. These extensions will be included in future work.
\end{rmk}

We also provide a characterization of the optimal policy for the infinite horizon problem described in Problem \ref{prob.IH}. This is stated as Proposition \ref{prop.IH} below.
\begin{prop}\label{prop.IH}
        Consider Problem \ref{prob.IH} 
for the system in \eqref{eq.FKE} described in Section \ref{sec.Model}, with stage cost \eqref{eq.stagecost}.
Then the policy $U^*$ in \eqref{eq.optsoIH} is an optimal policy.
        \begin{equation}\label{eq.optsoIH}
                U^* = \diag\left(\frac{u_{max}}{2}(\mathds{1}-\sgn(y^T{R_\mathcal{D}}+v^T))\right)
        \end{equation}
where $y^T\in \mathbb{R}^{n}$ is a solution of
\begin{align}\label{eq.paramStateIH}
&q^T + y^TR_\mathcal{E}\Gamma E  \\
&+ \frac{u_{max}}{2}\left((y^TR_\mathcal{D}+v^T)-\lvert y^TR_\mathcal{D}+v^T\rvert_{elem}\right)H=0\notag
\end{align}
\end{prop}
\begin{proof}
	See Section \ref{sec.Methods}.
\end{proof}
\begin{rmk}
	Proposition \ref{prop.IH} establishes that the optimal policy $U^*$ for Problem \ref{prob.IH} are constant and are independent of the state of network and of time.
\end{rmk}

\section{Examples}\label{sec.examples}
We compute the optimal policies for two 
examples using a numerical solution of differential equation \eqref{eq.paramState}. We have also validated our findings by comparing the optimal policies obtained from Propositions \ref{prop.FH} and \ref{prop.IH} with the ones obtained using the value iteration algorithm \cite{Bertsekas2005} for a discrete time approximation of \eqref{eq.FKE}. The time simulations presented have been performed using the Stochastic Simulation Algorithm \cite{Gillespie1977} implemented in the software package GillesPy2 \cite{Abel2016}.

\begin{figure}[h!]
        \centering
	\tikzset{
  queuei/.pic={
    \draw[line width=0.8pt]
      (0,0) -- ++(2cm,0) -- ++(0,-1cm) -- ++(-2cm,0);
    \foreach \Val in {1,...,4}
      \draw ([xshift=-\Val*8pt]2cm,0) -- ++(0,-1cm);
    \node[above] at (-0.50cm,-0.50cm) {$X_{#1}$};
    \node at (0.2cm,-0.5cm) {\Large $\cdot$};
    \node at (0.5cm,-0.5cm) {\Large $\cdot$};
    \node at (0.8cm,-0.5cm) {\Large $\cdot$};
  },
  mytri/.style={
    draw,
    shape=isosceles triangle,
    isosceles triangle apex angle=60,
    inner xsep=6pt
  }
}

\begin{tikzpicture}[>=latex]

 \path
    (0,0cm) pic {queuei=1}
    (0,-2cm) pic {queuei=2};

\draw (2.5,-0.5cm) circle [radius=0.5cm];
\draw (2.5,-2.5cm) circle [radius=0.5cm];

\draw[->] (3,-0.5) -- +(20pt,0);
\draw[->] (3,-2.5) -- +(20pt,0);
\node at (2.5,-0.5cm) {$\gamma_1 x_1$};
\node at (2.5,-2.5cm) {$\gamma_2 x_2$};

  \draw[->,bend right=25] (0.8,-1cm) to node[midway,right] {$u_1 x_1$} (0.8,-2cm);
\end{tikzpicture}
        \caption{Example network with $n=2$ queues. The associated routing matrices and state transition matrices are$R_{\mathcal{E}} =
        \begin{bmatrix}
                -1&0\\
                0&-1
\end{bmatrix}$,
$R_{\mathcal{D}} =
        \begin{bmatrix}
                -1\\
                 +1
        \end{bmatrix}$,
        $H =
        \begin{bmatrix}
                1&0
        \end{bmatrix}$,
        $E =
        \begin{bmatrix}
                1&0\\
                0&1
        \end{bmatrix}$
}\label{fig.2dchain}
\end{figure}
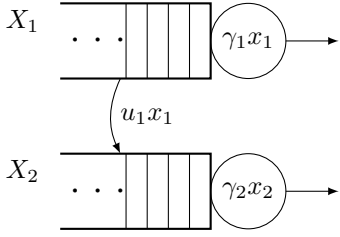

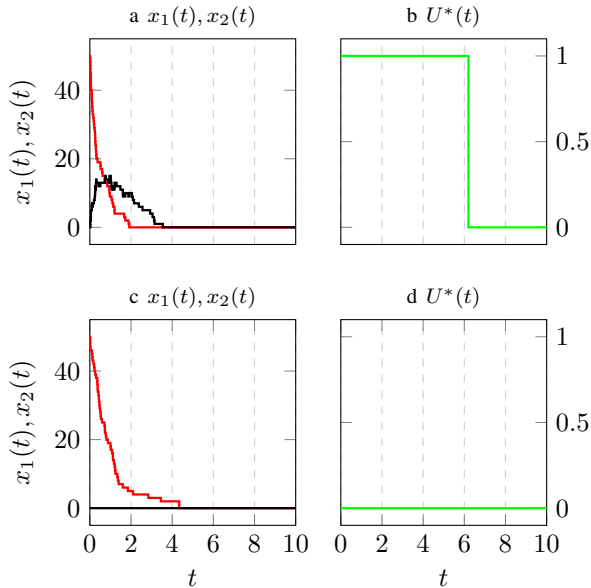
\begin{figure}[h!]
	\begin{tikzpicture}
\begin{groupplot}[
    group style             = {group size = 2 by 2, vertical sep = 1cm, horizontal sep = 0.6cm, xlabels at = edge bottom , ylabels at = edge left },
    view                    = {0}{90},
    width		    = 0.5\columnwidth,
    height		    = 0.5\columnwidth,
    ylabel 	            = {$x_1(t), x_2(t)$},
    xlabel 	            = {$t$},
    xmin                    = 0.0,
    xmax                    = 10,
    ]
\nextgroupplot[
		xmajorgrids,
		grid style 				= {thin, dashed, black!20},
                xmajorticks=false, 
		]
\addplot
		[
            		const plot,
			draw			= red,
			mark			= none,
			line width		= 0.03cm,
            style          = solid
		]
		table
		[
			x				= time,
			y				= state1,
		]
        {TrajectoryFirstExample.dat};\label{plt.FirstExState}
		\addplot
		[
			const plot,
			draw			= black,
			mark			= none,
			line width		= 0.03cm,
		        style          = solid
		]
		table
		[
			x				= time,
			y				= state2,
		]
        {TrajectoryFirstExample.dat};\label{plt.FirstExStatetwo}
\nextgroupplot[
                xmajorticks=false, 
		xmajorgrids,
		grid style 				= {thin, dashed, black!20},
		yticklabel pos=right
		]
\addplot
		[
            		const plot,
			draw			= green,
			mark			= none,
			line width		= 0.03cm,
            style          = solid
		]
		table
		[
			x				= time,
			y				= control,
		]
        {FirstExample.dat};\label{plt.FirstExControl}
\nextgroupplot[
		xmajorgrids,
		grid style 				= {thin, dashed, black!20},
		]
\addplot
		[
            		const plot,
			draw			= red,
			mark			= none,
			line width		= 0.03cm,
            style          = solid
		]
		table
		[
			x				= time,
			y				= state1,
		]
        {TrajectoryFirstExample2.dat};\label{plt.FirstExState2}
		\addplot
		[
			const plot,
			draw			= black,
			mark			= none,
			line width		= 0.03cm,
		        style          = solid
		]
		table
		[
			x				= time,
			y				= state2,
		]
        {TrajectoryFirstExample2.dat};\label{plt.FirstExStatetwo2}
\nextgroupplot[
		xmajorgrids,
		grid style 				= {thin, dashed, black!20},
		yticklabel pos=right,
		ymin = -0.1,
		ymax = 1.1
		]
	\addplot		[
            		const plot,
			draw			= green,
			mark			= none,
			line width		= 0.03cm,
            style          = solid
		] coordinates {(0,0) (10,0)};

\end{groupplot}
\node[text width=6cm, align=center, anchor = south] at (group c1r1.north) {\subcaption{$x_1(t), x_2(t)$\label{fig_optStateFirst}}};
\node[text width=6cm, align=center, anchor = south] at (group c2r1.north) {\subcaption{$U^*(t)$\label{fig_optControlFirst}}};
\node[text width=6cm, align=center, anchor = south] at (group c1r2.north) {\subcaption{$x_1(t), x_2(t)$\label{fig_optStateFirst2}}};
\node[text width=6cm, align=center, anchor = south] at (group c2r2.north) {\subcaption{$U^*(t)$\label{fig_optControlFirst2}}};
\end{tikzpicture}
\caption{The optimal policies $U^*$ (\ref{plt.FirstExControl}) for two different cases of cost coefficients are shown in the right hand side diagrams for the network in Fig.\ref{fig.2dchain}. The policies are bang-bang with respect to time validating Proposition \ref{prop.FH}. The left hand side diagrams show the time evolution of the states $x_1$ (\ref{plt.FirstExState}) and $x_2$ (\ref{plt.FirstExStatetwo}) from the initial condition $x=[50,0]^T$ in a particular trajectory. The diagrams in the top row are obtained for state cost $q^T=[2.5,1]$ and control costs $v=1$. The bottom row diagrams are obtained for $q^T=[2.5,1]$ and costs $v=2$.
}\label{fig.FirstExample}
\end{figure}

We first consider the network of 2 queues (\emph{i.e.} $n=2$) shown in Fig.\ref{fig.2dchain} with one routing control signal $u_1$ taking values from $0$ to $u_{max}=1$ that transition units from queue $X_1$ to queue $X_2$. The exit rate parameters of queue $X_1$ and $X_2$ are $\gamma_1=\gamma_2=1$. In Fig.\ref{fig.FirstExample} we make use of Proposition \ref{prop.FH} and solve Problem \ref{prob.FH} over a horizon $T=10$ with initial condition $x_0=[50,0]^T$.
Specifically, in Fig.\ref{fig_optControlFirst} and Fig. \ref{fig_optStateFirst} we compute the optimal policy and simulate the state evolution for $v=1$ and $q^T=[2.5,1]$. For these costs it is optimal to route units from the first queue to the second because of the high $x_1$ state cost and low control costs. The optimal policy takes the  value of $u_{max}=1$ until the queues have emptied. 

In Fig.\ref{fig_optControlFirst2} we compute for the same system the optimal policy for $v=2$ and in Fig.\ref{fig_optStateFirst2} we show the associated evolution of the state. The optimal policy in this case is to empty queue $X_1$ without using any routing due to the higher values of control cost, \emph{i.e.} $U^*(t) =0 \forall t$.
In Section \ref{sec.Appendix} we provide also a representative figure of the value function for the previous example.

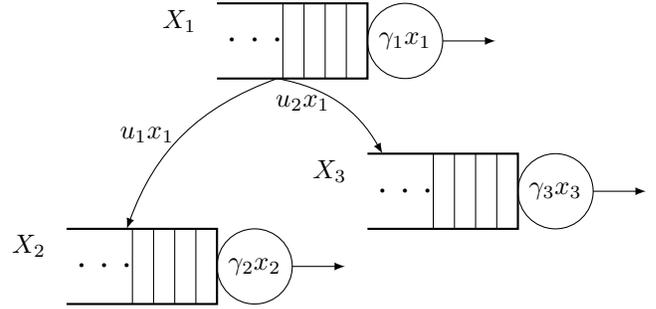
\begin{figure}[h!]
        \centering
\tikzset{
  queuei/.pic={
    \draw[line width=0.8pt]
      (0,0) -- ++(2cm,0) -- ++(0,-1cm) -- ++(-2cm,0);
    \foreach \Val in {1,...,4}
      \draw ([xshift=-\Val*8pt]2cm,0) -- ++(0,-1cm);
    \node[above] at (-0.50cm,-0.50cm) {$X_{#1}$};
    \node at (0.2cm,-0.5cm) {\Large $\cdot$};
    \node at (0.5cm,-0.5cm) {\Large $\cdot$};
    \node at (0.8cm,-0.5cm) {\Large $\cdot$};
  },
  mytri/.style={
    draw,
    shape=isosceles triangle,
    isosceles triangle apex angle=60,
    inner xsep=6pt
  }
}

\begin{tikzpicture}[>=latex]
 \path
    (0,0cm) pic {queuei=1}
    (-2,-3cm) pic {queuei=2}
    (2,-2cm) pic {queuei=3};

\draw (0.5,-3.5cm) circle [radius=0.5cm];
\draw (2.5,-0.5cm) circle [radius=0.5cm];
\draw (4.5,-2.5cm) circle [radius=0.5cm];

\draw[->] (3,-0.5) -- +(20pt,0);
\draw[->] (1,-3.5) -- +(20pt,0);
\draw[->] (5,-2.5) -- +(20pt,0);
\node at (2.5,-0.5cm) {$\gamma_1 x_1$};
\node at (0.5,-3.5cm) {$\gamma_2 x_2$};
\node at (4.5,-2.5cm) {$\gamma_3 x_3$};

  \draw[->,bend right=25] (0.8,-1cm) to node[midway,left] {$u_1 x_1$} (-1.2,-3cm);
  \draw[->,bend left=25] (0.8,-1cm) to node[midway,left] {$u_2 x_1$} (+2.2,-2cm);
\end{tikzpicture}
	\caption{Network with $n=3$ queues. The routing matrices and state transition matrices are $R_{\mathcal{E}} =
        \begin{bmatrix}
                -1&0&0\\
                0&-1&0\\
                0&0&-1
\end{bmatrix}$,
$R_{\mathcal{D}} =
        \begin{bmatrix}
                -1&0\\
                +1&-1\\
                 0&1
        \end{bmatrix}$,
        $H =
        \begin{bmatrix}
                1&0&0\\
                1&0&0
        \end{bmatrix}$,
        $E =
        \begin{bmatrix}
                1&0&0\\
                0&1&0\\
                0&0&1
        \end{bmatrix}$
}\label{fig.3dscheme}
\end{figure}

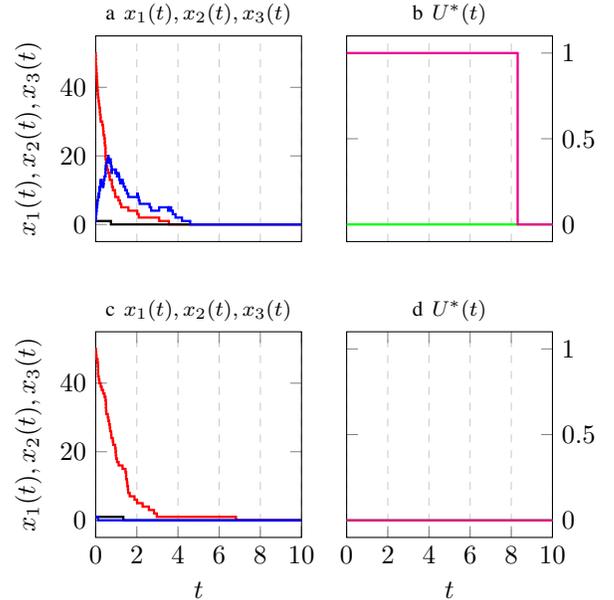
\begin{figure}[h!]
\begin{tikzpicture}
\begin{groupplot}[
    group style             = {group size = 2 by 2, vertical sep = 1.2cm, horizontal sep = 0.6cm, xlabels at = edge bottom , ylabels at = edge left },
    view                    = {0}{90},
    width		    = 0.5\columnwidth,
    height		    = 0.5\columnwidth,
    ylabel 	            = {$x_1(t), x_2(t), x_3(t)$},
    xlabel 	            = {$t$},
    xmin                    = 0.0,
    xmax                    = 10,
    ]
\nextgroupplot[ 
		xmajorgrids,
		grid style 				= {thin, dashed, black!20},
                xmajorticks=false, 
		]
\addplot
		[
            		const plot,
			draw			= red,
			mark			= none,
			line width		= 0.03cm,
            style          = solid
		]
		table
		[
			x				= time,
			y				= state1,
		]
        {TrajectorySecondExample.dat};\label{plt.SecondExState}
		\addplot
		[
			const plot,
			draw			= black,
			mark			= none,
			line width		= 0.03cm,
		        style          = solid 
		]
		table
		[
			x				= time,
			y				= state2,
		]
        {TrajectorySecondExample.dat};\label{plt.SecondExStatetwo}
		\addplot
		[
			const plot,
			draw			= blue,
			mark			= none,
			line width		= 0.03cm,
		        style          = solid 
		]
		table
		[
			x				= time,
			y				= state3,
		]
        {TrajectorySecondExample.dat};\label{plt.SecondExStatethree}
\nextgroupplot[ 
                xmajorticks=false, 
		xmajorgrids,
		grid style 				= {thin, dashed, black!20},
		yticklabel pos=right
		]
\addplot
		[
            		const plot,
			draw			= green,
			mark			= none,
			line width		= 0.03cm,
            style          = solid
		] coordinates {(0,0) (10,0)};\label{plt.secondexu1}
\addplot
		[
            		const plot,
			draw			= magenta,
			mark			= none,
			line width		= 0.03cm,
            style          = solid
	    ] coordinates {(0,1) (8.3,1) (8.31,0) (10,0)};\label{plt.secondexu2}
\nextgroupplot[ 
		xmajorgrids,
		grid style 				= {thin, dashed, black!20},
		]
\addplot
		[
            		const plot,
			draw			= red,
			mark			= none,
			line width		= 0.03cm,
            style          = solid
		]
		table
		[
			x				= time,
			y				= state1,
		]
        {TrajectorySecondExampleNocontrol.dat};
		\addplot
		[
			const plot,
			draw			= black,
			mark			= none,
			line width		= 0.03cm,
		        style          = solid 
		]
		table
		[
			x				= time,
			y				= state2,
		]
        {TrajectorySecondExampleNocontrol.dat};
		\addplot
		[
			const plot,
			draw			= blue,
			mark			= none,
			line width		= 0.03cm,
		        style          = solid 
		]
		table
		[
			x				= time,
			y				= state3,
		]
        {TrajectorySecondExampleNocontrol.dat};
\nextgroupplot[ 
		xmajorgrids,
		grid style 				= {thin, dashed, black!20},
		yticklabel pos=right,
		ymin = -0.1,
		ymax = 1.1
		]
\addplot[
            		const plot,
			draw			= green,
			mark			= none,
			line width		= 0.03cm,
            style          = solid
		] coordinates {(0,0) (10,0)};
	\addplot		[
            		const plot,
			draw			= magenta,
			mark			= none,
			line width		= 0.03cm,
            style          = solid
		] coordinates {(0,0) (10,0)};

\end{groupplot}
\node[text width=6cm, align=center, anchor = south] at (group c1r1.north) {\subcaption{$x_1(t), x_2(t),x_3(t)$\label{fig_optStateSecond}}}; 
\node[text width=6cm, align=center, anchor = south] at (group c2r1.north) {\subcaption{$U^*(t)$\label{fig_optControlSecond}}};
\node[text width=6cm, align=center, anchor = south] at (group c1r2.north) {\subcaption{$x_1(t), x_2(t),x_3(t)$\label{fig_optStateSecond2}}}; 
\node[text width=6cm, align=center, anchor = south] at (group c2r2.north) {\subcaption{$U^*(t)$\label{fig_optControlSecond2}}};
\end{tikzpicture}
\caption{We show the optimal policies $U^*=\diag(u_1^*,u_2^*)$ ($u_1$\ref{plt.secondexu1}, $u_2$\ref{plt.secondexu2}) for two costs configuration in the right hand side diagrams for the network in Fig.\ref{fig.3dscheme}. The policies are bang-bang with respect to time validating Proposition \ref{prop.FH} and show that it is optimal to route units from queue one to queue three. The left hand side diagrams show the time evolution of the states $x_1(t)$ (\ref{plt.SecondExState}), $x_2(t)$(\ref{plt.SecondExStatetwo}) and $x_3(t)$(\ref{plt.SecondExStatethree}) from the initial condition $x=[50,1,1]^T$. The diagrams in the top row are obtained for state cost coefficient $q^T=[2.5,1,1]$ and control cost coefficient $v=1.6$. The bottom row diagrams are obtained for $q^T=[2.5,1,1]$ and costs $v=2.2$.}\label{fig.SecondEx}
\end{figure}

In Fig. \ref{fig.SecondEx} we compute the optimal policies and state evolution for Problem \ref{prob.FH} using the result provided in Proposition \ref{prop.FH} for a network with 3 queues (\emph{i.e.} $n=3$) shown in Fig. \ref{fig.3dscheme}.
Also for this system the optimal policies are bang-bang with respect to time as suggested by Proposition \ref{prop.FH}. The diagrams in Fig.\ref{fig_optStateSecond} and Fig.\ref{fig_optControlSecond} are obtained for state cost coefficient $q^T=[2.5,1,1]$ and control cost coefficient $v=1.6$. The diagrams in the lower row are obtained for $v=2$. For low values of control cost $v$ it is optimal to route units from queue $X_1$ to queue $X_3$. It should be noted that also in this case it is optimal to route units to queue $X_3$ until the system reaches equilibrium.

\section{Derivation of the results}\label{sec.Methods}
We start by considering Problem \ref{prob.FH} and we make use of 
dynamic programming principles 
\cite{Bellman2013} to characterise the optimal policies. We then give the proofs of Proposition~\ref{prop.FH} and Proposition \ref{prop.IH} in Section \ref{s.proofFH} and Section \ref{s.proofIH}, respectively.

In order to derive our results we make use of the \acf{HJBE} a partial difference equation for the value function $\mathcal{V}(x,t):\mathbb{Z}^n_{\geq}\times \mathbb{R}_\geq \to \mathbb{R}_\geq$ of an optimal control problem, which provides a sufficient condition for optimality \cite{Hanson2007}.

In the derivations we denote $\mathcal{U}^d$ the set of diagonal matrices with elements in $\mathcal{U}$, \emph{i.e.} $\mathcal{U}^d:=\{\text{diag}(u_1, \hdots, u_{m_u}): u_k\in\mathcal{U} \  \forall k\}$.

The \ac{HJBE} for Problem \ref{prob.FH} associated with the \ac{MJP} in \eqref{eq.FKE} as described in Section \ref{sec.Model} is \cite{Theodorou2012}
\begin{align}\label{eq.HJBE}
	\min_{U
	\in\mathcal{U}^d}&
	\left[\dfrac{\partial \mathcal{V}(x,t)}{\partial t}+g_c(x(t),U)\right.\\
	&\left.+\sum_{i\in \mathcal{I}}W_i(x,t)(\mathcal{V}(x+r_i,t)-\mathcal{V}(x,t)) \right]=0 \notag
\end{align}
subject to the boundary condition
\begin{equation}\label{eq.HJBEBoundary}
	\mathcal{V}(x(T),T) = \tilde g_c(x(T))
\end{equation}

We also define the \ac{HJBE} for Problem \ref{prob.IH} 
\begin{align}\label{eq.HJBEIH}
	\min_{U\in\mathcal{U}^d}
&\left[g_c(x(t),U)\right.\\
	&\left.	+\sum_{i\in \mathcal{I}}W_i(x,t)(\mathcal{V}(x+r_i,t)-\mathcal{V}(x,t)) \right]=0
\end{align}

We give the proof for our main result in Proposition \ref{prop.FH} in Section \ref{s.proofFH} and the proof for Proposition \ref{prop.IH} in Section \ref{s.proofIH}.
\subsection{Proof of Proposition \ref{prop.FH}}\label{s.proofFH}
\begin{proof}
Consider \eqref{eq.HJBE} and let us rearrange it as
\begin{align}\label{eq.HJBERef0}
	\min_{U\in\mathcal{U}^d}&\left[\dfrac{\partial \mathcal{V}(x,t)}{\partial t}+g_c(x(t),U)\right.\\
 &\left.+\sum_{i\in \mathcal{D}}W_i(x,t)(\mathcal{V}(x+r_i,t)-\mathcal{V}(x,t))\right.\notag\\
 &\left.+\sum_{i\in \mathcal{E}}W_i(x,t)(\mathcal{V}(x+r_i,t)-\mathcal{V}(x,t))\right]=0\notag
\end{align}
We now define the following quantities
\begin{align}\label{eq.vectorFiniteVE}
&\Delta_{R_\mathcal{E}}\mathcal{V}(x,t)= [\mathcal{V}(x+r_k,t)-\mathcal{V}(x,t)]_{k\in \mathcal{E}} \\
&\Delta_{R_\mathcal{D}}\mathcal{V}(x,t)=[\mathcal{V}(x+r_k,t)-\mathcal{V}(x,t)]_{k\in \mathcal{D}\label{eq.vectorFiniteVD}}
\end{align}
that are row vectors where each element is the finite difference of the value function defined with a corresponding column $r_k$ of the state change matrices $R_\mathcal{D}, R_\mathcal{E}$ (these are associated with the events in $\mathcal{D}, \mathcal{E}$ we consider).

We make use of the  quantities defined in \eqref{eq.vectorFiniteVE}, \eqref{eq.vectorFiniteVD}, the definition of matrices $\Gamma$, $H$, $E$ and $U$ and by substituting the expression for $g_c$ in \eqref{eq.stagecost}, \eqref{eq.HJBERef0} can be written as
\begin{align}\label{eq.HJBERef1}
& \dfrac{\partial \mathcal{V}(x,t)}{\partial t}+ q^Tx
+\Delta_{R_\mathcal{E}}\mathcal{V}(x,t) \Gamma E x \\
&+ \min_{U\in\mathcal{U}^d}\left[(\Delta_{R_\mathcal{D}}\mathcal{V}(x,t)+v^T)U\right]Hx=0\notag
\end{align}
We now explicitly compute the minimisation appearing in the previous expression yielding the optimal $U^*(x,t)$
\begin{align}\label{eq.minop}
	U^*(x,t)&=\argmin_{U\in\mathcal{U}^d}\left[(\Delta_{R_\mathcal{D}}\mathcal{V}(x,t)+v^T)U\right]\\
&= \diag\left(\frac{u_{max}}{2}(\mathds{1}-\sgn(\Delta_{R_\mathcal{D}}\mathcal{V}(x,t)+v^T))\right)\notag
\end{align}
where the second equality holds by making use of the definition of $\sgn$ and of matrix $U$. 
\nomenclature{\ $\lvert x\rvert_{elem}$}{Element wise absolute of $v$: $\lvert v\rvert_{elem}= [\lvert v_1 \rvert,\ldots,  \lvert v_n \rvert]$}

By making use of the previous equation we obtain the following nonlinear partial difference equation 
\begin{align}\label{eq.HJBERef3}
&\dfrac{\partial \mathcal{V}(x,t)}{\partial t}+q^Tx
+ \Delta_{R_\mathcal{E}}\mathcal{V}(x,t) \Gamma E x \\\notag
&+\frac{u_{max}}{2}\left[(\Delta_{R_\mathcal{D}}\mathcal{V}(x,t)+v^T)\right. \\
&\left.-|\Delta_{R_\mathcal{D}}\mathcal{V}^T(x,t)+v^T|_{elem}\right]Hx =0 \notag
\end{align}
We now consider the following candidate value function as solution for \eqref{eq.HJBERef3}
\begin{equation}\label{eq.candidateV}
	\mathcal{V}(x,t) = y^T(t)x
\end{equation}
where $y\in \mathbb{R}^{n} \quad \forall t \in \mathbb{R}_\geq$.
We now substitute \eqref{eq.candidateV} in \eqref{eq.vectorFiniteVE} and \eqref{eq.vectorFiniteVD} yielding
\begin{align}
&\Delta_{R_\mathcal{E}}\mathcal{V}(x,t)= y^T(t)R_\mathcal{E}
\label{eq.vectorFiniteVEMat}
\\
&
\Delta_{R_\mathcal{D}}\mathcal{V}(x,t)=y^T(t)R_\mathcal{D}
\label{eq.vectorFiniteVDMat}
\end{align}
We now substitute \eqref{eq.candidateV} and the expressions above in \eqref{eq.HJBERef3} and obtain
\begin{align}\label{eq.HJBERef4}
&\left[\dot y^T(t)+q^T + y^T(t)R_\mathcal{E}\Gamma E+\right. \\
& \left. \frac{u_{max}}{2}\left((y^T(t)R_\mathcal{D}+v^T)-\lvert y^T(t)R_\mathcal{D}+v^T\rvert_{elem}\right)H\right]x = 0\notag
\end{align}
This holds for all values of $x$ when the differential equation in \eqref{eq.paramState} holds.
Note that we also have the terminal condition $y(T)=c$ from \eqref{eq.finalcost}, \eqref{eq.HJBEBoundary}, \eqref{eq.candidateV} as stated in Proposition \ref{prop.FH}.

Therefore the optimal policy $U^*(x,t)$ is obtained by using \eqref{eq.minop} and \eqref{eq.candidateV} yielding
\begin{equation}
U^*(x,t)= \diag\left(\frac{u_{max}}{2}(\mathds{1}-\sgn(y^T(t){R_\mathcal{D}}+v^T))\right)
\end{equation}
where $y^T(t)$ satisfies \eqref{eq.HJBERef4} with the terminal condition $y(T)=c$. 
\end{proof}

\subsection{Proof of Proposition \ref{prop.IH}}\label{s.proofIH}
\begin{proof}
The \ac{HJBE} in \eqref{eq.HJBEIH} takes the form 
\begin{align}\label{eq.HJBEIHRef0}
	\min_{U\in\mathcal{U}^d}&\left[q^Tx + v^TUHx\right.\\
 &\left.+\sum_{i\in \mathcal{D}}W_i(x,t)(\mathcal{V}(x+r_i)-\mathcal{V}(x))\right.\notag\\
&\left.+\sum_{i\in \mathcal{E}}W_i(x,t)(\mathcal{V}(x+r_i)-\mathcal{V}(x))\right]=0\notag
\end{align}
We re-write \eqref{eq.HJBEIHRef0} as follows
\begin{align}\label{eq.HJBERef1IH}
	q^Tx+&\Delta_{R_\mathcal{E}}\mathcal{V}(x) \Gamma E x \\
&+ \min_{U\in\mathcal{U}^d}\left[(\Delta_{R_\mathcal{D}}\mathcal{V}(x)+v^T)U\right]Hx=0\notag
\end{align}
where we have made use of the following quantities
\begin{align}\label{eq.vectorFiniteVEIH}
&\Delta_{R_\mathcal{E}}\mathcal{V}(x)= [\mathcal{V}(x+r_k)-\mathcal{V}(x)]_{k\in \mathcal{E}} \\
&\Delta_{R_\mathcal{D}}\mathcal{V}(x)=[\mathcal{V}(x+r_k)-\mathcal{V}(x)]_{k\in \mathcal{D}\label{eq.vectorFiniteVDIH}}
\end{align}
We now compute the minimisation appearing in \eqref{eq.HJBERef1IH} yielding
\begin{align}\label{eq.minopIH}
	U^*&=\argmin_{U\in\mathcal{U}^d}\left[(\Delta_{R_\mathcal{D}}\mathcal{V}(x)+v^T)U\right]\\
&= \diag\left(\frac{u_{max}}{2}(\mathds{1}-\sgn(\Delta_{R_\mathcal{D}}\mathcal{V}(x)+v^T))\right)\notag
\end{align}
and by substituting this expression back in \eqref{eq.HJBERef1IH} after some manipulation we obtain
\begin{align}\label{eq.HJBEIHRef3}
&q^Tx
+ \Delta_{R_\mathcal{E}}\mathcal{V}(x) \Gamma E x \\\notag
&+\frac{u_{max}}{2}\left[(\Delta_{R_\mathcal{D}}\mathcal{V}(x)+v^T)\right. \\
&\left.-\lvert\Delta_{R_\mathcal{D}}\mathcal{V}(x)+v^T\rvert_{elem}\right]Hx =0
\end{align}
We now make the following ansatz for the value function
\begin{equation}\label{eq.CandidateIH}
	\mathcal{V}(x) = y^Tx
\end{equation}
We make use of the previous expression in \eqref{eq.HJBEIHRef3} and we obtain
\begin{align}\label{eq.HJBERef4IH}
&\left[q^T +y^TR_\mathcal{E}\Gamma E\right. \\
& +\left. \frac{u_{max}}{2}\left((y^TR_\mathcal{D}+v^T)-\lvert y^TR_\mathcal{D}+v^T\rvert_{elem}\right)H\right]x = 0\notag
\end{align}
which holds for all values of $x$ when \eqref{eq.paramStateIH} holds, thus verifying the ansatz in \eqref{eq.CandidateIH}.

By substituting \eqref{eq.HJBERef4IH} in \eqref{eq.minopIH} we find the optimal policy is
\begin{align}\label{eq.minopIH}
	U^*&= \diag\left(\frac{u_{max}}{2}(\mathds{1}-\sgn(y^TR_\mathcal{D}+v^T))\right)\notag
\end{align}
where $y^T$ satisfies the non-linear equation \eqref{eq.HJBERef4IH}.
\end{proof}
\section{Conclusions and future work}\label{sec.Conclusions}

We have considered the problem of finding optimal policies in stochastic networks of $M/M/\infty$ queues, with controlled routing events and exit events with prescribe rates.
We have explicitely characterised the optimal policies for costs that are linear functions of the system states and the control inputs. We have also shown that the optimal policies do not involve state feedback, whereby the control input is adjusted based on the state of the network, but they depend on the network topology and the system parameters. Computations of optimal policies validate our findings. Future work includes incorporating larger class of events and constraints.

\section{Appendix}\label{sec.Appendix}
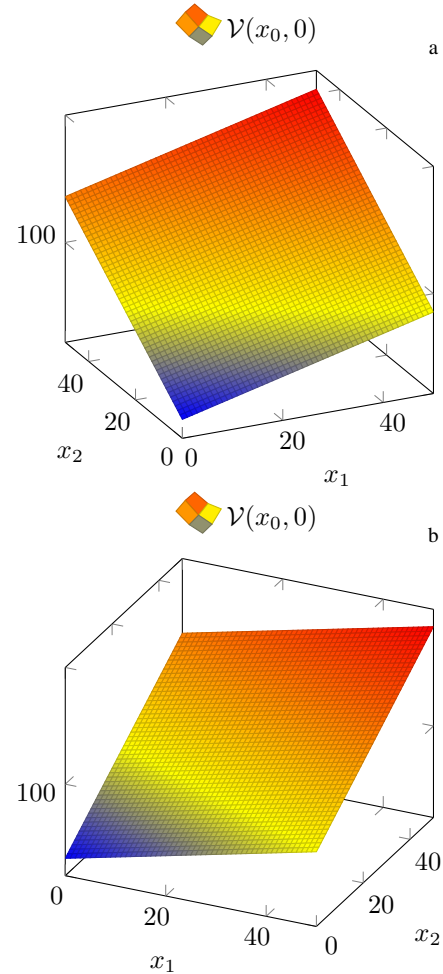
\begin{figure}[!h]
\begin{center}
\begin{tikzpicture}
\begin{groupplot}[
    group style             = {group size = 1 by 2, vertical sep = 1.6cm, horizontal sep = 0cm },
    width		    = 0.75\columnwidth,
    height		    = 0.75\columnwidth,
    ]
\nextgroupplot[
		legend style={at={(0.5,+1.05)},anchor=south,draw=none},
    		ylabel 	            = {$x_2$},
		xlabel 	            = {$x_1$},
    		view                    = {335}{25},
		]
\addplot3[
	surf,
	]
	table{ValueFunction2d.dat};
	\addlegendentry{$\mathcal{V}(x_0,0)$};
\nextgroupplot[
	legend style={at={(0.5,+1.05)},anchor=south,draw=none},
    		xlabel 	            = {$x_1$},
    		ylabel 	            = {$x_2$},
		]
\addplot3[
	surf,
	]
	table{ValueFunction2d.dat};
	\addlegendentry{$\mathcal{V}(x_0,0)$};
\end{groupplot}
\node[text width=0.2cm, align=center, anchor = south] at (group c1r1.north east) {\subcaption{\label{VF1}}};
\node[text width=0.2cm, align=center, anchor = south] at (group c1r2.north east) {\subcaption{\label{VF2}}};
\end{tikzpicture}
\end{center}
	\caption{The diagrams in Fig.\ref{VF1} and Fig.\ref{VF2} show the value function $\mathcal{V}(x_0,0)$ for the example considered in Fig.\ref{fig.2dchain} from two different azimuth and elevation angles. The value function is a linear function of the states that is inline the findings discussed in Section \ref{sec.Methods}.}\label{fig.2dValueFunction}
\end{figure}

In Fig. \ref{fig.2dValueFunction} we show the optimal value function $\mathcal{V}(x_0,0)$ for the example network in Fig. \ref{fig.2dchain} discussed in Section \ref{sec.examples}. The function is linear in the variable $x$ as discussed in the derivation of our results (see \eqref{eq.candidateV}). The optimal value function $\mathcal{V}$ in Fig. \ref{fig.2dValueFunction} was obtained by numerical integration of \eqref{eq.paramState}. The result of the numerical integration has also been validated by comparing it to the solution obtained using the value iteration algorithm for a discrete time approximation of the \ac{HJBE} \eqref{eq.HJBERef4}. 


\begin{thebibliography}{10}
\providecommand{\url}[1]{#1}
\csname url@samestyle\endcsname
\providecommand{\newblock}{\relax}
\providecommand{\bibinfo}[2]{#2}
\providecommand{\BIBentrySTDinterwordspacing}{\spaceskip=0pt\relax}
\providecommand{\BIBentryALTinterwordstretchfactor}{4}
\providecommand{\BIBentryALTinterwordspacing}{\spaceskip=\fontdimen2\font plus
\BIBentryALTinterwordstretchfactor\fontdimen3\font minus
  \fontdimen4\font\relax}
\providecommand{\BIBforeignlanguage}[2]{{%
\expandafter\ifx\csname l@#1\endcsname\relax
\typeout{** WARNING: IEEEtran.bst: No hyphenation pattern has been}%
\typeout{** loaded for the language `#1'. Using the pattern for}%
\typeout{** the default language instead.}%
\else
\language=\csname l@#1\endcsname
\fi
#2}}
\providecommand{\BIBdecl}{\relax}
\BIBdecl

\bibitem{Erlang1909}
A.~K. Erlang, ``The theory of probabilities and telephone conversations,''
  \emph{Nyt Tidsskrift for Matematik B}, 1909.

\bibitem{Crabill1972}
T.~B. Crabill, ``Optimal control of a service facility with variable
  exponential service times and constant arrival rate,'' \emph{Management
  Science}, vol.~18, no.~9, pp. 560--566, May 1972.

\bibitem{Stidham1985}
S.~Stidham, ``Optimal control of admission to a quenching system,'' \emph{IEEE
  Transactions on Automatic Control}, vol.~30, no.~8, pp. 705--713, Aug. 1985.

\bibitem{Johansen1980}
S.~G. Johansen and S.~Stidham, ``Control of arrivals to a stochastic
  input–output system,'' \emph{Advances in Applied Probability}, vol.~12,
  no.~4, pp. 972--999, Dec. 1980.

\bibitem{Topkis1978}
D.~M. Topkis, ``Minimizing a submodular function on a lattice,''
  \emph{Operations Research}, vol.~26, no.~2, pp. 305--321, Apr. 1978.

\bibitem{Li2019}
Q.-L. Li, J.-Y. Ma, R.-N. Fan, and L.~Xia, ``An overview for markov decision
  processes in queues and networks,'' 2019.

\bibitem{arazi2004}
A.~Arazi, E.~Ben-Jacob, and U.~Yechiali, ``Bridging genetic networks and
  queueing theory,'' \emph{Physica A: Statistical Mechanics and its
  Applications}, vol. 332, pp. 585--616, 2004.

\bibitem{Jackson1957}
J.~R. Jackson, ``Networks of waiting lines,'' \emph{Operations Research},
  vol.~5, no.~4, pp. 518--521, Aug. 1957.

\bibitem{Ghoneim1985}
H.~A. Ghoneim and S.~Stidham, ``Control of arrivals to two queues in series,''
  \emph{European Journal of Operational Research}, vol.~21, no.~3, pp.
  399--409, Sep. 1985.

\bibitem{Rosberg1982}
Z.~Rosberg, P.~Varaiya, and J.~Walrand, ``Optimal control of service in tandem
  queues,'' \emph{IEEE Transactions on Automatic Control}, vol.~27, no.~3, pp.
  600--610, Jun. 1982.

\bibitem{Hajek1984}
B.~Hajek and R.~G. Ogier, ``Optimal dynamic routing in communication networks
  with continuous traffic,'' \emph{Networks}, vol.~14, no.~3, pp. 457--487,
  Sep. 1984.

\bibitem{BadianPessot2019}
P.~Badian-Pessot, M.~E. Lewis, and D.~G. Down, ``Optimal control policies for
  an m/m/1 queue with removable server and dyanmic serivice rates,''
  \emph{Probability in the Engineering and Informational Sciences}, vol.~35,
  no.~2, pp. 189--209, Jul. 2019.

\bibitem{Gallisch1979}
E.~Gallisch, ``On monotone optimal policies in a queueing model ofm/g/1 type
  with controllable service time distribution,'' \emph{Advances in Applied
  Probability}, vol.~11, no.~4, pp. 870--887, Dec. 1979.

\bibitem{HernandezLerma1983}
O.~Hernández-Lerma and S.~I. Marcus, ``Adaptive control of service in queueing
  systems,'' \emph{Systems \& Control Letters}, vol.~3, no.~5, pp. 283--289,
  Nov. 1983.

\bibitem{Weber1987}
R.~R. Weber and S.~Stidham, ``Optimal control of service rates in networks of
  queues,'' \emph{Advances in Applied Probability}, vol.~19, no.~1, pp.
  202--218, Mar. 1987.

\bibitem{Su2024}
Y.~Su and J.~Li, ``Admission control of double-sided queues with multiple
  customer types,'' \emph{IEEE Transactions on Automatic Control}, vol.~69,
  no.~3, pp. 1960--1966, Mar. 2024.

\bibitem{Dimitrakopoulos2017}
Y.~Dimitrakopoulos and A.~Burnetas, ``The value of service rate flexibility in
  an m/m/1 queue with admission control,'' \emph{IISE Transactions}, vol.~49,
  no.~6, pp. 603--621, Apr. 2017.

\bibitem{Ata2006}
B.~Ata and S.~Shneorson, ``Dynamic control of an m/m/1 service system with
  adjustable arrival and service rates,'' \emph{Management Science}, vol.~52,
  no.~11, pp. 1778--1791, Nov. 2006.

\bibitem{Adusumilli2010}
K.~M. Adusumilli and J.~J. Hasenbein, ``Dynamic admission and service rate
  control of a queue,'' \emph{Queueing Systems}, vol.~66, no.~2, pp. 131--154,
  Sep. 2010.

\bibitem{Feinberg2015}
E.~A. Feinberg and X.~Zhang, ``Optimal switching on and off the entire service
  capacity of a parallel queue,'' \emph{Probability in the Engineering and
  Informational Sciences}, vol.~29, no.~4, pp. 483--506, Oct. 2015.

\bibitem{Li2024}
Y.~Li and A.~Rantzer, ``Exact dynamic programming for positive systems with
  linear optimal cost,'' \emph{IEEE Transactions on Automatic Control},
  vol.~69, no.~12, pp. 8738--8750, Dec. 2024.

\bibitem{Blanchini2023}
F.~Blanchini, P.~Bolzern, P.~Colaneri, G.~De~Nicolao, and G.~Giordano,
  ``Optimal control of compartmental models: The exact solution,''
  \emph{Automatica}, vol. 147, p. 110680, Jan. 2023.

\bibitem{Blackwell1964}
D.~Blackwell, ``Memoryless strategies in finite-stage dynamic programming,''
  \emph{The Annals of Mathematical Statistics}, vol.~35, no.~2, pp. 863--865,
  jun 1964.

\bibitem{Puterman1994}
M.~L. Puterman, Ed., \emph{Markov Decision Processes}.\hskip 1em plus 0.5em
  minus 0.4em\relax John Wiley {\&} Sons, Inc., apr 1994.

\bibitem{Bertsekas2005}
D.~Bertsekas, \emph{Dynamic programming and optimal control}.\hskip 1em plus
  0.5em minus 0.4em\relax Belmont, Mass: Athena Scientific, 2005.

\bibitem{Gillespie1977}
D.~T. Gillespie, ``Exact stochastic simulation of coupled chemical reactions,''
  \emph{The Journal of Physical Chemistry}, vol.~81, no.~25, pp. 2340--2361,
  Dec. 1977.

\bibitem{Abel2016}
J.~H. Abel, B.~Drawert, A.~Hellander, and L.~R. Petzold, ``Gillespy: A python
  package for stochastic model building and simulation,'' \emph{IEEE Life
  Sciences Letters}, vol.~2, no.~3, pp. 35--38, Sep. 2016.

\bibitem{Bellman2013}
R.~Bellman, \emph{Dynamic Programming}, ser. Dover Books on Computer
  Science.\hskip 1em plus 0.5em minus 0.4em\relax Newburyport: Dover
  Publications, 2013, description based upon print version of record.

\bibitem{Hanson2007}
F.~B. Hanson, \emph{Applied stochastic processes and control for
  Jump-diffusions}, ser. Advances in design and control.\hskip 1em plus 0.5em
  minus 0.4em\relax Philadelphia, Pa.: Society for Industrial and Applied
  Mathematics (SIAM, 3600 Market Street, Floor 6, Philadelphia, PA 19104),
  2007, no.~13, includes bibliographical references (p. 403-422) and index. -
  Description based on title page of print version.

\bibitem{Theodorou2012}
E.~A. Theodorou and E.~Todorov, ``Stochastic optimal control for nonlinear
  markov jump diffusion processes,'' in \emph{2012 American Control Conference
  (ACC)}.\hskip 1em plus 0.5em minus 0.4em\relax IEEE, Jun. 2012, pp.
  1633--1639.

\end{thebibliography}

\end{document}